\documentclass[12pt]{amsart}

\usepackage{geometry}
\usepackage{amssymb}

\usepackage{graphicx}
\usepackage [english]{babel}
\usepackage [autostyle, english = american]{csquotes}
\MakeOuterQuote{"}

\makeatletter
\@namedef{subjclassname@2020}{
 \textup{2020} Mathematics Subject Classification}
\makeatother

\newtheorem{theorem}{Theorem}[section]
\newtheorem{corollary}[theorem]{Corollary}
\newtheorem{lemma}[theorem]{Lemma}
\newtheorem{proposition}[theorem]{Proposition}

\theoremstyle{definition}
\newtheorem{definition}[theorem]{Definition}

\numberwithin{equation}{section}

\newcommand{\df}[1]{\ensuremath{\mathrm{Def}}(#1)}
\newcommand{\Def}[1]{\ensuremath{\mathrm{Def^n}}(#1)}

\newenvironment{claim}
{\begin{trivlist}  \item \textsc{Claim}~} {\end{trivlist}}
\newenvironment{proof of claim}
{\begin{trivlist}  \item \textsc{Proof of Claim:}~} {\hfill $\Box$ (\textsc{Claim})
\end{trivlist}}

\def \S{\textbf{S}}
\def \M{\textbf{M}}
\newtheorem{qu}[theorem]{Question}

\def \Bn {\operatorname{B[n]}}

\title{%Approximate cancellation in Grothendieck groups of first-order structures
A definable matching theorem and approximate Grothendieck semigroups}

\author{Jana Ma\v{r}\'{i}kov\'{a}}
\address{Kurt G\"{o}del Research Center for Mathematical Logic, Universit\"{a}t Wien, 1090 Wien,
Austria}

\begin{document}

\begin{abstract}
For an arbitrary first-order structure, we introduce the notion of an approximate Grothendieck semigroup, which satisfies a Banach-Schr\"oder-Bernstein Theorem, as well as a weak Cancellation 
Law (wCL). 

At the heart of the proof of the wCL lies a definable matching theorem, which vastly generalizes a related result in Ma\v{r}\'{i}kov\'{a} \cite{mat}, 
and is of independent interest.
\end{abstract}


%\keywords{aaaa, bbbb, cccc}

\maketitle

\maketitle

\begin{section}{Introduction}
 Throughout this paper, $\S=(S;\dots)$ denotes a first-order structure with universe $S$.  {\em Definable\/} shall mean definable in $\S$, with parameters. We denote by $\Def{\S}$ the boolean algebra of definable subsets of $S^n$.

We will be concerned with an instance of the following question. Let $G$ be a group of bijections $S^n \to S^n$, and $\mathcal{B}$ a $G$-invariant boolean algebra of subsets of $S^n$ (i.e. if $g\in G$ and $X\in \mathcal{B}$, then $g(X)\in \mathcal{B}$).  Is there a finitely additive $G$-invariant measure on $\mathcal{B}$ that assigns 1 to a distinguished element from $\mathcal{B}$? 

There are, of course, many well-studied instances of the above.  For example, when $\mathcal{S}=(S;1,\cdot )$ is any group, $G:=\{ \tau_s \colon s\in S \mbox{ and } \tau_s (x)=s\cdot x \mbox{ for all }x\in S  \}$, and $\mathcal{B}$ is the power set of $S$, with the distinguished set being $S$ itself, then this is just asking whether $\S$ is amenable as a discrete group.  

In model theory, one is interested in situations where there are some definability requirements on $G$ and/or $\mathcal{B}$. 
For instance, $G$ might be an $\mathcal{S}$-definable group acting on the boolean algebra of its $\mathcal{S}$-definable subsets by left translation, with the distinguished set being $G$. Then $G$ is called definably amenable if the answer to the above question is yes.

A non-amenable group might well be definably amenable, as for example the free group on two generators is well-known to be non-amenable, however, it is definably amenable, since its theory is stable (Sela \cite{sela}), and stable groups are definably amenable (Poizat \cite{poizat}).

 Here, we consider the case when $G$ is a group of definable bijections $S^n \to S^n$ (so $G$ is a subgroup of the group of all definable bijections $S^n \to S^n$), whose elements we shall refer to as {\em isomorphisms\/}, and we set $\mathcal{B}=\Def{\S}$.  The requirement that $G$ is a group of total definable bijections is largely cosmetic. Our results would for instance go through in the more general setting, where $G$ is a pseudogroup of definable bijections (see Ceccherini-Silberstein et al. \cite{cgh} for an account of amenability in the setting of actions of pseudogroups).

The motivation for this work stems from an open question in o-minimality.  O-minimal structures are certain densely linearly ordered structures, which are particularly well-behaved from a toplogical point of view, with the semialgebraic sets serving as an archetypical example. Another example is the structure $\mathbb{R}_{\textup{an}, \exp}$, the real field with restricted analytic functions and the exponential, whose o-minimality has famously played an important part in applications to diophantine geometry and Hodge theory in work of J. Pila, J. Tsimerman, and many others.

While much of differential calculus goes through in any o-minimal expansion of a real closed field (definable functions are piecewise differentiable, derivatives of definable functions are definable, there is a definable version of the Implicit Function Theorem etc.), the question of integration is more subtle.  In general, an o-minimal structure is not closed under taking antiderivatives, however, in some cases they are known to be definable in an o-minimal expansion of the original structure (see for instance Lion, Rolin \cite{LionRolin} for the globally subanalytic case, or Kaiser \cite{K1} for a theory of integration for non-archimedean real closed fields). 
But even the following, more basic question, a variation of which was posed by Hrushovski, Peterzil, Pillay \cite{nip}, has yet to be answered in full generality.

\noindent
\begin{qu}\label{qu_omin}
Let $\S$ be an o-minimal expansion of a real closed field, and let $G$ consist of the definable $C^1$-diffeomorphisms $S^n \to S^n$ with Jacobian determinant equal to $\pm 1$ everywhere.  Let $X \in \Def{\S}$ be bounded, with non-empty interior.  Is there a $G$-invariant finitely additive measure on the boolean algebra of bounded definable subsets of $S^n$ which takes values in $\mathbb{R}^{\geq 0} \cup \{ \infty \}$ and assigns a positive value to $X$? 
\end{qu}

The answer is known to be positive in certain cases, for instance when $\S$ satisfies the common theory of all $\mathcal{L}$-expansions of the real field, where $\mathcal{L}$ is the language of $\S$, in which case the theory of $\S$ is called {\em pseudo-real\/}.  However, not all o-minimal expansions of real closed fields are pseudo-real - see Lipschitz, Robinson \cite{LR} and Hurshovski, Peterzil \cite{hp}.  

We note that finitely additive measures on boolean algebras of definable sets are called {\em Keisler measures\/}. They first appeared in Keisler \cite{keisler}, and rose to prominence after their use in the proof of the o-minimal group conjecture \cite{nip}.  Keisler measures play an important part especially in NIP theories (i.e. in theories where every definable family has finite VC-dimension), in some sense taking over the role of types in stable theories - see work by Chernikov, Hrushovski, Pillay, Simon, Starchenko, and many others.

One way of addressing questions regarding the existence of invariant Keisler measures, is by
using a connection between the existence of finitely additive measures and the non-existence of so-called paradoxical decompositions that was established by Tarski \cite{tarski1},\cite{tarski}.
Below, the relation $\leq $ denotes the usual pre-order on a commutative monoid: $\alpha \leq \beta$ if there is $\gamma$ such that $\alpha + \gamma = \beta$.  We use $n \alpha$ as an abbreviation for the expression $\alpha + \dots +\alpha$, in which $\alpha$ appears $n$ times. 
\begin{theorem}[Tarski]\label{tarski} Let $(C; 0,+,\alpha )$ be a commutative monoid with identity element 0 and a distinguished element $\alpha$.  Then the following are equivalent.
\begin{enumerate}
\item For every $k$, $(k+1)\alpha \not\leq k\alpha$. 
\item There is a homomorphism of monoids $\mu \colon C \to [0,\infty]$ such that $\mu \alpha =1$. 
\end{enumerate}
\end{theorem}
Returning to Question \ref{qu_omin}, we set \[\Bn:= \{Y \in \Def{\S} \colon Y \mbox{ is bounded}\},\] and take $C$ to be the quotient monoid $\Bn/_{\simeq }$, where $+$ on $\Bn$ is disjoint union, $\simeq$ is the equivalence relation induced by ismomorphisms, and $\alpha := [X]_{\simeq}$. Then the homomorphism $\mu$ in $(2)$ is the desired invariant Keisler measure on $\Bn$. 

Note that Condition $(1)$ can be phrased as a statement about the non-existence of a definable matching covering the first part of 
a definable bipartite graph with parts $(k+1)X$ and $kX$, and an edge relation given by finitely many isomorphisms $g_1 , \dots ,g_k$, where $x\in (k+1)X$, $y\in kX$ are adjacent if $g_i (x)=y$ for some $i\in \{1,\dots ,k\}$.  We shall refer to such graphs as {\em mp graphs\/} (see Definition \ref{mp}).
Tarski's theorem can be applied in a similar way outside of the bounded o-minimal setting, though one needs to ensure the existence of disjoint unions.

So the question arises, under what conditions mp graphs admit definable matchings covering one part of the vertex set. 
We remark that the Hall-Rado-Hall infinite Marriage Theorem (see e.g. \cite{wagon}, p. 323) fails in the definable setting, and not even the stronger condition of regularity guarantees the existence of a perfect definable matching, by a well-known example of Laczkovich \cite{L}.

In \cite{mat}, an approximate definable matching theorem has been proved for 2-regular mp graphs in the special case where $\mathcal{S}$ is an o-minimal expansion of a real closed field.  While the strategy in \cite{mat} follows ideas from Elek, Lippner \cite{el} and Lyons, Nazarov \cite{ln}, the proofs rely heavily on tools from o-minimality.  Moreover, useful information is only obtained for definable sets $X$ whose standard part has nonempty interior.

Here, we prove an approximate matching theorem (Theorem \ref{matchingthm1}) which holds in any first-order structure: 
Let $k\geq 2$.  Then, in a $k$-regular mp graph, for every $m>0$, we can find a definable matching that covers the vertex set $V$ outside of a subset that is $m$-small in $V$.  A set is $m$-small in $V$ if, roughly, $m$ disjoint copies of it embed isomorphically into $V$.
The strategy follows again \cite{ln}. That is, we first show the existence of definable matchings without short augmenting paths (this part of the proof works without assuming regularity), and  
subsequently, we use regularity to derive a bound on the size of the set of vertices not covered by $M$.  In the second part of the argument it is convenient to use Tarski's Theorem \ref{tarski}: If one part of the vertex set satisfies the negation of the first condition, then any definable matching covers that part up to a subset, which is $m$-small in it, for every $m>0$.  If it does satisfy the first condition, then we have an invariant Keisler measure at hand.

Using isomorphisms to compare sizes of sets leads naturally to the notion of an approximate Grothendieck semiring, which is defined in Section 5. 
There, we establish a Banach-Schr\"oder-Bernstein Theorem (BSB), where BSB refers to a generalization by Banach of the usual Schr\"oder-Bernstein Theorem of set theory.
Furthermore, using a version of Theorem \ref{matchingthm1} for certain definable multigraphs, we establish a weak Cancellation Law (wCL) for the approximate Grothendieck semigroup.
The remaining two sections discuss some consequences for definable versions of equidecomposability and Tarski's theorem.

For a reference on model theory, see for instance Hodges \cite{hodges}; a possible reference for graph theory is Diestel \cite{diestel}.

\smallskip\noindent
{\bf Acknowledgements.\/} We thank K. Peterzil for discussions, in particular for pointing out that some preliminary proofs apply in a much wider setting.  We thank I. Kaplan for enlightening questions.  We thank J. Kraj\'i\v{c}ek for pointing us to his work on aEC.

This work was partially supported by the FWF Principal Investigator Project (DOI 10.55776/PAT4567425).
\end{section}

\begin{section}{Notation and definitions} 
A function $f\colon A\to S^n$, $A\subseteq S^m$, is said to be definable, if its graph \[\Gamma f = \{ (x,y) \in S^{m+n} \colon f(x)=y \}\] is definable.

A {\em graph\/} is a pair of sets $(V,E)$, where $V\not=\emptyset$ is called the set of vertices, and $E\subseteq V^2$ is a symmetric and antireflexive relation whose elements are called edges. So graphs are not oriented, have no loops, and we shall think of edges, when convenient, as two-element subsets of $V$.  In this sense, a vertex $v$ and an edge $e$ are {\em incident\/} if $v\in e$. Two vertices $v,w$ are {\em adjacent\/} if $(v,w)\in E$. The degree of a vertex is the number of edges incident with it.  If $v\in e$, then we also say that $e$ {\em covers\/} $v$. 
If $Y\subseteq E$, then $Y(V)$ denotes the set of vertices covered by the edges in $Y$. 
A graph is {\em definable\/} if $V$ and $E$ are definable. 
A set $Y\subseteq E$ is {\em independent\/} if no two distinct edges in $Y$ share a vertex. A graph is called {\em $k$-regular\/}, if each vertex is incident with exactly $k$ edges.

A graph $(V,E)$ is {\em bipartite\/} if $V=A\dot\cup B$ and $E\subseteq A\times B$, where the pair $A$, $B$ is called a {\em bipartition\/} of the graph.  A {\em definable bipartite graph\/} is a definable graph with a definable bipartition.

 A {\em matching\/} in a graph $(V,E)$ is an independent set of edges.  
 A matching is {\em perfect\/} if it covers $V$. Given a matching $M\subseteq E$, a path in $(V,E)$ is {\em alternating\/} if its edges alternate belonging to $M$. An {\em augmenting path\/} is a 
 simple
 alternating path such that the starting vertex and the final vertex are not covered by $M$. The length $|p|$ of a path $p$ is its number of edges, and so an augmenting path is necessarily of odd length.  
 A $k$-(augmenting) path is a (augmenting) path of length $k$.
 To {\em flip\/} an augmenting path $p$ means to remove the edges of $p$ that were in $M$ from $M$, and to place the edges of $p$ that were not in $M$, into $M$.  Note that, after flipping an augmenting path, the new $M$ is still a matching.

The letters $i,j,k,l,m,n,K,N$ are used to denote natural numbers, where $\mathbb{N}=\{ 0,1,2,\dots \}$.
We use $[k]$ to denote the set $\{ 1,\dots ,k \}$.
\end{section}

\begin{section}{Matchings without short augmenting paths}
{\bf Throughout this paper, we fix a group $G$ of definable bijections $S^n \to S^n$, and we refer to its elements as isomorphisms.\/} 
\begin{definition}\label{mp}
An {\em mp graph\/} is a definable bipartite graph with bipartition $A,B \subseteq S^n$ and edge relation $E\subseteq A\times B$, such that 
there is a definable finite partition $\{ A_1 ,\dots ,A_l \}$ of $A$, and for each $i\in [l]$, \[E \cap (A_i \times S^n ) = \dot\bigcup_{j\in \mathcal{F}_i} \Gamma f_{ij} , \] where $\mathcal{F}_i$ is finite, and each $f_{ij}$ is an isomorphism restricted to $A_i$. 
\end{definition}
Note that for every mp graph, there is a positive integer which is an upper-bound for the degrees of its vertices.
Also note that for every path $(x_0 , x_1 , \dots ,x_k )$ starting in $A$ there is a unique sequence $(f_1  , \dots ,f_k )$ of restricted isomorphisms from $\bigcup_{i\in [l]} \{f_{ij}\colon j \in \mathcal{F}_i \}$ such that $x_{2i+1} =f_{2i+1} (x_{2i})$ and $x_{2i+2} = f^{-1}_{2i+2} (x_{2i + 1} )$.  We shall call $(f_1 , \dots ,f_k )$ the {\em generating sequence\/} of the path $(x_0 , \dots, x_k )$.

\begin{proposition}\label{noshortaps}
Let $(V,E)$ be an mp graph with bipartition $A, B \subseteq S^n$.
Given a definable matching $M_0 \subseteq E$ and $k\geq 1$, there is a definable matching $M \subseteq E$, which does not admit any augmenting paths of length $\leq k$, and is such that $M_0 (V) \subseteq M(V)$.
\end{proposition}
\begin{proof}
Let $\{ A_1 , \dots ,A_l \}$ and $\{ \mathcal{F}_1 , \dots ,\mathcal{F}_l \}$ be as in Definition \ref{mp}.  
The proof proceeds by induction on $k$.  At each stage, we show that, 
after flipping a certain set of augmenting paths, a) we obtain a matching, and  b) we did not introduce any new augmenting paths of equal or shorter length than the flipped ones.

\begin{trivlist}
  \item{{\em Base step.\/}}
  Suppose that for some $i\in [l]$, $A_i$ contains the starting vertex of a 1-augmenting path with respect to $M_0$, with generating sequence $f\in \mathcal{F}_i$. Let $X$ be the (definable) set of all vertices in $A_i$ that are starting vertices of 1-augmenting paths with generating sequence $f$.
  Then $M_0 \cup \Gamma f|_{X}$ is a definable matching, because $\Gamma f|_{X}$ only contains edges between vertices not covered by $M_0$, and each vertex in $X \cup f(X )$ is incident with exactly one edge in $\Gamma f$, as $f$ is a bijection. Furthermore, all vertices that were covered by $M_0$ are still covered, there are no more 1-augmenting paths with generating sequence $f$ starting in $A_i$, and we did not introduce any new 1-augmenting paths.
  
  Since $\mathcal{F}_i$ is finite, it suffices to use the above observation finitely many times to arrive at a matching that covers $V(M_0 )$, does not admit any 1-augmenting paths starting in $A_i$, and does not introduce any new 1-augmenting paths. We repeat this procedure for the remaining $A_j$'s to arrive at a definable matching covering $V(M_0)$ and not admitting any 1-augmenting paths. 

\item{{\em Inductive step.\/}}
  Now assume that $M_0 \subseteq E$ is a definable matching which does not admit any augmenting paths of length $\leq k-1$, where $k> 1$ is even.  We wish to construct a matching covering $V(M_0 )$ and not admitting any augmenting paths of length $\leq k+1$.

  Let $i$ be the first index in $[l]$ such that $A_i$ contains the starting vertex of a $(k+1)$-augmenting path, say with generating sequence $(f)=(f_1 ,f_2 ,\dots ,f_{k+1})$, and let $X$ be the (definable) subset of $A_i$ of starting vertices of $(k+1)$-augmenting paths with generating sequence $(f)$. We denote the set of these augmenting paths by $P_{X,(f)}\subseteq V^{k+2}$.  
  
\begin{claim}
Let $p,q \in P_{X,(f)}$ with $p\not=q$. Then $p$, $q$ do not share a vertex.
\end{claim}
\begin{proof of claim}
First observe that $p_0$ cannot appear in $q$, since it would have to be $q_0$, but $p\not=q$. But then $p_{k+1}$ cannot appear in $q$, since it would have to be $q_{k+1}$, but that is impossible because $p_0 \not=q_0$, and all functions appearing in $(f)$ are bijections. 
Now suppose $p_i = q_j$ for some $i,j \in [k]$. Again because $p_0 \not=q_0$ and all functions in $(f)$ are bijections, $i\not=j$. So we may assume $i<j$.  But then $(p_0 ,\dots ,p_{i-1}, q_{j}, \dots ,q_{k+1})$ is an augmenting path of length $<k+1$, a contradiction.
\end{proof of claim}
If follows that the relation $M'$ obtained from $M_0$ by flipping the paths in $P_{X,(f)}$ is a matching.

\begin{claim}
Suppose $p$ is an augmenting path in $(V,E)$ with respect to $M'$, of length $\leq k+1$. Then $p$ is an augmenting path for $M_0$.
\end{claim}
\begin{proof of claim}

If $p$ is a 1-augmenting path for $M'$, then $p$ was already an augmenting path with respect to $M_0$, since flipping augmenting paths cannot uncover vertices.
So we may assume that $|p|\geq 3$.  We assert towards a contradiction that $p$ is an augmenting path for $M'$ but not for $M_0$.  Then $p$ contains at least one edge that was flipped as part of a path in $\mathcal{P}_{X,(f)}$.  Let $c^1 , \dots ,c^l$ be the subpaths of $p$ of maximal length that were flipped as part of a path in $\mathcal{P}_{X,(f)}$, listed in the order of their occurrence in $p$, i.e. if $i<j$, then every edge in $c^i$ occurs in $p$ before every edge in $c^j$.  By maximality, $p$ contains at least one edge between any consecutive pair of components $c^i$, $c^{i+1}$. 
For each $i$, let $q^i$ be the unique path in $\mathcal{P}_{X,(f)}$ such that $c^i$ is a subpath of $q^i$.    
To simplify notation, we shall proceed assuming $l=3$, but the same proof works for any $l\geq 1$.  

We write $p$ as a concatenation of paths $a^1 c^1 a^2 c^2 a^3 c^3 a^4$, where the $a^j$'s are uniquely determined by the $c^i$'s. Furthermore, we write $q^i = q^{i1}c^i q^{i2}$ for $i=1,2,3$.

\includegraphics[height=2.5cm,width=12cm,angle=360]{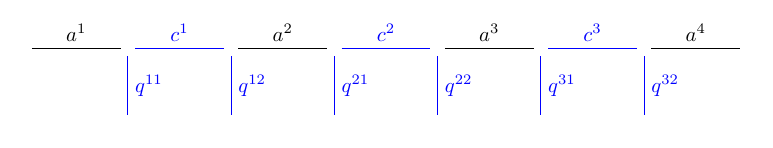}

Then 
\begin{equation}\label{1} \sum_{i=1}^{4}|a^i | + \sum_{i=1}^{3}|c^i| \leq k+1 \end{equation}
and 
\begin{equation}\label{2} |q^{i1}|+|c^i |+|q^{i2}| = k+1 \mbox{ for each } i \in \{ 1,2,3 \}.
\end{equation} 

We aim to show that then at least one of  
\[
a^1 q^{11}, \;
q^{12}a^{2}q^{21}, \;
q^{22}a^{3}q^{31}, \;
q^{32} a^4
\]
is an augmenting path for $M_0$ of length $\leq k-1$, which will yield a contradiction with the inductive assumption. Clearly, all these paths are augmenting paths with respect to $M_0$. Suppose to the contrary that the four augmenting paths above are all of length $>k-1$.  Then they are all of length $\geq k+1$, since they are of odd length.
So, by adding up their lengths, \[ \sum_{i=1}^{4}|a^i| + \sum_{i=1}^{3}\sum_{j=1}^{2}|q^{ij}| \geq 4(k+1),\]
and hence by equations \ref{1} and \ref{2},
\[ 4(k+1) - 2\sum_{i=1}^{3} |c_i| \geq 4(k+1),\]
thus
$ \sum_{i=1}^{3}c_i \leq 0$, contradicting the assumption that at least one edge in $p$ was flipped when flipping $\mathcal{P}_{X,(f)}$.
\end{proof of claim}

Since there are only finitely many generating sequences for $(k+1)$-augmenting paths starting in $A_i$, using the above two claims, we may assume that $M'$ is a definable matching which does not admit any $(k+1)$-augmenting paths starting in $A_i$, and which covers $V(M_0 )$. Each subsequent $A_{j}$ is handled similarly. After finitely many steps we thus arrive at a definable matching, which does not admit any augmenting paths of length $\leq k+1$ starting in $A$ (hence there are no such paths starting in $B$ either), and which covers $M_0 (V)$.
\end{trivlist}

\end{proof}

\end{section}

\begin{section}{Matchings in $k$-regular bipartite graphs}
We wish to construct definable matchings which are almost perfect (this is the best we can hope for), using isomorphisms to compare "sizes" of sets.  The main result here is Theorem \ref{matchingthm1}, of which we discuss several variants, to be used in subsequent sections.

We shall need to ensure the existence of disjoint unions.  To achieve this, we use a set-up from \cite{nip}: 
An $m$-cycle is a finite formal sum $\Sigma_{i\in [m]} X_i$ of $X_i \in \Def{\S}$.  In particular, each $X_i$ is a 1-cycle. 
Let $\textbf{M}$ be the free commutative monoid generated by $\Def{\S}$ with (formal) addition of cycles.

\begin{definition}
\begin{enumerate}
\item A {\em pw-embedding\/} is a map \[ f\colon \Sigma_{i\in [m]} X_i \to \Sigma_{j\in [n]} Y_j \] between two cycles such that for every $i\in [m]$, 
there is a finite definable partition $\{ X_{ik} \}$ of $X_i$ such that for every $k$, $f(X_{ik}) \subseteq Y_j$ for some $j\in [n]$, and $f|_{X_{ik}}$ is the restriction of an  isomorphism to $X_{ik}$, and $f$ is injective in the following sense: if $f(X_{ik})\subseteq Y_j$ and $f(X_{i'k'})\subseteq Y_j$ and $f(X_{ik})\cap f(X_{i'k'})\not=\emptyset$, then $i=i'$ and $k=k'$.
\item A {\em pw-isomorphism} is a pw-embedding which is onto. 
\item For two cycles $X=\Sigma_{i\in [k]}X_i$, $X'=\Sigma_{i\in [k]}X'_i$ we write $X\subseteq X'$, and call $X$ a {\em subcycle\/} of $X'$, if $X_i \subseteq X'_i$ for each $i\in [k]$, and $X\setminus X'$ is the cycle $\Sigma_{i\in [k]}(X_i \setminus X'_i )$.
\end{enumerate}
\end{definition}

\begin{definition}\label{ordering}
Let $X=\Sigma_{i\in [k]} X_i , Y=\Sigma_{j\in [l]}Y_j \in \M$.
\begin{enumerate}
\item We write $X \simeq Y$ if there is a pw-isomorphism 
$f\colon X \to Y$. If additionally $X,Y \in \Def{\bf{S}}$, then we say that $X,Y$ are {\em equidecomposable\/}. 

\item We write $X \leq Y$ if there is a pw-embedding 
$f\colon X \to Y$.

\item By $nX$ we mean the cycle \[\Sigma_{i\in [k]} X_i + \dots + \Sigma_{i\in [k]}X_i \] where each $X_i$ appears $n$ times.

\item For $m>0$, we write $X \leq_m Y$ if there are $p,q \in \mathbb{N}^{>0}$ such that $m \leq \frac{p}{q}$ and $p X \leq q Y$. In this case we also say that $X$ is {\em $m$-small\/} in $Y$.

\item We write $X \leq_0 Y$ if for each $m>0$, there is a cycle $X' \subseteq X$ such that $X' \leq_m X$ and $X \setminus X'  \leq Y$. 

\item We write $X \simeq_0 Y$ if for each $m>0$ 
there are cycles $X' \subseteq X$ and $Y' \subseteq Y$ such that $X' \leq_m X$, $Y' \leq_m Y$ and $X \setminus X' \simeq Y \setminus Y'$.

\end{enumerate}
\end{definition}
The relations $\leq$ and $\leq_0$ are pre-orders on $\M$.  In the next section we shall see that $\leq_0$ is in fact an order.
Also note that if there is an isomorphism-invariant Keisler measure $\nu$ on $\mathcal{S}$, and $X,Y$ are $1$-cycles, then $ X\leq_m Y$ implies $m \cdot \nu X \leq \nu Y$.  
%Also note that for $X, X_1 , \dots ,X_k \in \M$, if $X_1 \leq_{km} X , \dots ,X_k \leq_{km} X$, then $X_1 + \dots + X_k \leq_m X$.

\smallskip\noindent
The next theorem says that every $k$-regular $mp$ graph contains a definable matching covering the vertex set outside of an arbitrarily small set. 

In the proof of the theorem below, for $X\subseteq A \dot\cup B$ and $M\subseteq E$, we set
\[N_M (X):= \{ y\in A \dot\cup B \colon \, \exists x\in X \ \; (x,y)\in M  \} . \]

\begin{theorem}\label{matchingthm1}
Let $(A\dot\cup B,E)$ be a $k$-regular mp graph, where $k\geq 2$.  
Then, for every $m>0$, there is a definable matching $M\subseteq E$ covering the vertex set outside of a definable set $Y_{0A}\dot\cup Y_{0B}$, where $Y_{0A}\subseteq A$, $Y_{0B}\subseteq B$ and $Y_{0A} \leq_{m} A$, $Y_{0B}\leq_{m} B$.
\end{theorem}
\begin{proof}
Fix an odd $K>0$.
By Prop \ref{noshortaps},
we can find a definable matching $M\subseteq E$ which does not admit any augmenting paths of length $\leq K$. Let $Y_{0A}\dot\cup Y_{0B}$ be the set of vertices not covered by $M$, where $Y_{0A} \subseteq A$ and $Y_{0B}\subseteq B$. Our goal is to 
show that for sufficiently large $K$, the matching $M$ is as required.  We shall compute the bound only for $Y_{0A}$, as the proof for $Y_{0B}$ is analogous.

Throughout the remainder of the proof, we let $i$ range over $\{ 0, \dots ,K+1 \}$. We define $Y_i$ as follows. We set $Y_0 := Y_{0A}$.  If $i$ is even, then $Y_{i+1} =N_G (Y_{i})$.   If $i$ is odd, then $Y_{i+1}=N_{M} (Y_{i})$.

Note that if there is $k>0$ such that $(k+1)A \leq k A$, then any definable matching $M\subseteq E$ works, since then $A\leq_m A$ for every $m>0$.  Else, we can set $C:=\M/_{\simeq}$ and $\alpha := A/_{\simeq}$ in 
Theorem \ref{tarski}, and thus obtain a homomorphism of monoids $\mu \colon \M/_{\simeq} \to [0,\infty ]$ such that $\mu A=1$.

Every vertex of $Y_{0}$ is incident with $k$ vertices from $Y_{1}$, and $Y_{1}$ is covered by $M$, since there are no 1-augmenting paths. It follows that   
every vertex in $Y_{1}$ is incident with at most $k-1$ many vertices from $Y_{0}$.
 
So 
\begin{equation}\label{0and1}
k Y_{0} \leq (k-1)Y_{1} .
\end{equation}
For $i$ odd, $Y_{i}$ is covered by $M$, by the non-existence of augmenting paths of length $\leq i$. Hence
$Y_{i} \simeq Y_{i+1 }$.

Note that $Y_3 = Y_1 + Y'_3$, where $Y'_{3}$ consists of the vertices in $Y_3$ that are not adjacent via an edge in $M$ to a vertex in $Y_2$. By the non-existence of short augmenting paths, every vertex in $Y'_3$ is covered by $M$, so can be adjacent to at most $(k-1)$ vertices from $Y_2$. Since every vertex in $Y_1$ is adjacent to a vertex in $Y_{0}$, it can again be adjacent to only at most $(k-1)$ vertices in $Y_2$.  So we obtain
\[ kY_{1} = k Y_{2} \leq (k-1)Y_3 = (k-1)(Y_1 + Y'_3 ), \mbox{ hence } \]
\[\mu Y_1 \leq (k-1) \mu Y'_{3} .\]
It is easy to see that in general, for $1<i\leq K$ odd, we have
\[ Y_i = Y_1 + Y'_3 + \dots + Y'_i ,\] where $Y'_i$ consists of the vertices in $Y_i$ that are not adjacent via an edge in $M$ with a vertex in $Y_{i-1}$.  And, by similar reasoning as in the case of $Y_3$,
\[ k(Y_1 + Y'_3 + \dots + Y'_{i-1} ) \leq (k-1)(Y_1 + Y'_i ) + k(Y'_3 + \dots + Y'_{i-1}).\]  It follows that 
\begin{equation}\label{1toidoubleprime}
\mu Y_1 \leq (k-1) \mu Y'_{i},
\end{equation}
hence $\frac{1}{k-1} \mu Y_1 \leq \mu Y'_i$ for every odd $i$. So
\[ 1 \geq \mu Y_{K+1}= \mu Y_K = \mu (Y_1 + Y'_3 + \dots + Y'_K ) \geq \mu Y_1 + \frac{K-1}{2} \cdot \frac{1}{k-1} \cdot \mu Y_1 . \] Then
\[ \mu Y_1 \leq \frac{2(k-1)}{2(k-1)+K-1}, \mbox{ and so }\]
\[ \mu Y_{0} \leq \frac{k-1}{k} \mu Y_1 \leq \frac{2(k-1)^2 }{k(2k+K-3)}. \]

\end{proof}
We shall need a one-sided approximate version of the above theorem for mp multigraphs (to be defined below).  We first prove it for mp graphs.
\begin{corollary}\label{matchingthm2}
Let $2\leq k \leq l$, and let $(A\dot\cup B,E)$ be an mp graph, which is of degree $\geq k$ in $A\setminus A'$, where $A' \subseteq A$ is such that $A' \leq_{l} A$, and of degree $\leq k$ in $B$.  
Then, for every positive rational $q<\frac{l}{k}$, there is a definable matching $M\subseteq E$ covering $A \setminus Y_{0}$, where $Y_{0} \leq_{q} A$.
\end{corollary}
\begin{proof}
The proof proceeds just as the proof of Theorem \ref{matchingthm1}, and so we only point out how to adjust the inequalities.  Inequality \ref{0and1} becomes
\[ kY_{0}  \leq (k-1) Y_1 + k A',\]
and Inequality \ref{1toidoubleprime} turns into
\[ \mu Y_1 \leq (k-1)\mu Y'_{i} + k \mu A' .\]
Then 
\[ 1\geq \mu Y_K = \mu (Y_1 + Y'_3 + \dots + Y'_K ) \geq \mu Y_1 + \frac{K-1}{2} \cdot (\frac{1}{k-1} \cdot \mu Y_1 - \frac{k}{k-1}\cdot \mu A' ), \mbox{ hence}\]
\[ Y_1 \leq \frac{k(K-1)}{K+2k-3} \cdot \frac{1}{l},\] 
where the right-hand side converges to $\frac{k}{l}$ from above as $K\to \infty$, and so 
\[ Y_{0} \leq \frac{(k-1)}{k} \frac{k(K-1)}{K+2k-3} \cdot \frac{1}{l} + \frac{1}{l}, \]
where, again, the quantity on the right approaches $\frac{k}{l}$ from above as $K\to \infty$.
\end{proof}

Corollary \ref{matchingthm2} can be generalized to a statement about mp {\em multigraphs\/}, which will appear in the proof of the Weak Cancellation Law. These are defined in the same way as definable mp graphs, except that one is additionally given $N>0$ and a definable symmetric map $I\colon E \to \{ 1,\dots ,N \}$.  For $e\in E$, the value $I(e)$ is called the {\em multiplicity\/} of $e$. The degree of a vertex is then the number of edges incident with it, counting multiplicities.  So an edge $(v,w)$ with $I(v,w)=i$ contributes $i$ towards the degree of $v$.  Matchings and augmenting paths are defined in the same way as for graphs.
\begin{corollary}\label{multigraph} Let $2\leq k \leq l$, and let
    $(A\dot\cup B,E,I)$ be an mp multigraph which is of degree $\geq k$ in $A\setminus A'$, where $A'\subseteq A$ is such that $A' \leq_l A$, and each vertex in $B$ is of degree $\leq k$. Then, for every positive integer $q< \frac{l}{k}$ there is a definable matching $M\subseteq E$ covering $A$ outside of a subset $Y_0$ with $Y_{0} \leq_q A$.
    \end{corollary}
\begin{proof}
Let $K>0$.  We can apply Proposition \ref{noshortaps} to the mp graph $(A\dot\cup B,E)$, to  obtain a definable matching $M$ without augmenting paths of length $\leq K$.  By identifying an edge $(u,v) \in M$ with exactly one of the edges between $u$ and $v$ in the original multigraph, $M$ becomes an idependent set of edges in the multigraph.  It is also easy to check that $M$ does not admit any augmenting paths of length $\leq K$ in the multigraph. The computation of the upper bound on the size of the set of vertices not covered by $M$ then works just as in the proof of Corollary \ref{matchingthm2}. 
\end{proof}
\end{section}

\begin{section}{Banach-Schr\"oder-Bernstein and Weak Cancellation}
Our first aim is to show that $\M$ satisfies an approximate BSB theorem.
Given cycles $X=\Sigma_{i\in [k]}X_i$ and $Y=\Sigma_{i\in [k]}Y_i$, we write $X\cap Y$ for the cycle $\Sigma_{i\in [k]} (X_i \cap Y_i)$. The cycles $X,Y$ are said to be disjoint if $X\cap Y=\emptyset$. 
\begin{theorem}[BSB]\label{aBSB}
Suppose $A\leq_0 B$ and $B \leq_0 A$ for $A,B \in \M$.  Then $A \simeq_0 B$.
\end{theorem}
\begin{proof}
Fix $m>0$.  We wish to find cycles $V\subseteq A$ and $W \subseteq B$ with $V\leq_m A$ and $W\leq_m B$, and a pw isomorphism $\phi \colon A\setminus V \to B\setminus W$.

Set $k=3m$. By the assumption, there are cycles $A_1 , A', B_1 , B'$ with $A'\leq_k A$ and $B'\leq_k B$ and pw isomorphisms $f\colon A\setminus A' \to B\setminus B_1$ and $g\colon A\setminus A_1 \to B\setminus B'$.

We shall define a sequence $V_1 , V_2 ,\dots ,V_k$ of pairwise disjoint subcycles of $A$, and $W_0 , W_1 , \dots ,W_k$ a sequence of pairwise disjoint subcycles of $B$. To simplify notation, we set \[ V'_i := V_i \cap A' \mbox{ and } V''_i :=V_i \setminus A',  \] and similarly \[ W'_i := W_i \cap B' \mbox{ and } W''_i :=W_i \setminus B'.  \]
We set $W_0 := B_1 \setminus B'$ and $V_1 := g^{-1}(W_0)$.  Suppose $W_0 , \dots ,W_{i-1}$ and $V_1 , \dots ,V_i$ have already been defined.  Then $W_i := f(V''_i )$ and $V_{i+1}:= g^{-1} (W''_{i} )$.
Below, we shall abuse notation and write cycles as sums of cycles which are not necessarily 1-cycles.
\begin{claim} We have
$W_0 \simeq \Sigma_{i\in [k]} V'_i + \Sigma_{i\in [k]} W'_i + W''_k$, where $\Sigma_{i\in [k]} W'_i$ is $k$-small in $B$ and $W''_k$ is $(k+1)$-small in $B$, and the $V'_i$'s, as well as the $W'_i$'s are pairwise disjoint.   
\end{claim}
\begin{proof of claim}
 Equidecomposability and disjointness of the unions follows easily from the construction.  The cycle $\Sigma_{i\in [k]} W'_i$ is $k$-small in $B$ because it is a subcycle of $B'$, which is $k$-small in $B$. The pairwise disjoint cycles $W_0 , \dots ,W_k \subseteq B$ contain a copy of $W''_k$ each, so $W''_k$ is $(k+1)$-small in $B$.  
\end{proof of claim}
We now define a pw-isomorphism \[\phi \colon (A\setminus A') + \Sigma_{i\in [k]} V'_i \to (B\setminus B_1) + \Sigma_{i\in [k]} h^{(i-1)} g (V'_i),\] where $h^{(i-1)}$ is $h=gf^{-1}$ composed $(i-1)$ times with itself: 
Set $\phi|_{A\setminus A'} = f$, and for $i \in [k]$ we set $\phi|_{V'_i} = h^{(i-1)} g$.  Then $\phi$ and \[V := A'\setminus \Sigma_{i=1}^{k}V'_i \mbox{,  } W:= B_1 \setminus \Sigma_{i=1}^{k} \phi (V'_i)\] are as required: Since $V$ is a subcycle of $A'$, it is $m$-small in $A$. Moreover, it is easy to check that $\phi (V'_1 ), \dots ,\phi (V'_k )$ are pairwise disjoint subcycles of $W_0$.   We have 
\[
\begin{array}{rll}
W&=& (W_0 + (B_1 \cap B')) \setminus \Sigma_{i\in [k]} \phi (V'_i)\\
&=& (W_0 \setminus \Sigma_{i \in [k]} \phi (V'_i )) +(B_1 \cap B'),
\end{array}
\]
where $B_1 \cap B'$ is $k$-small in $B$, and, by the construction, \[W_0 \setminus \Sigma_{i \in [k]} \phi (V'_i )\simeq \Sigma_{i\in [k]}W'_i + W''_k ,\] 
with $\Sigma_{i\in [k]} W'_i$, as well as $W''_k$, $k$-small in $B$ by the claim. So by our choice of $k$, $W$ is $m$-small in $B$.
\end{proof}
Note that the above shows that $\leq_0$ is an order on $\M$.

\smallskip\noindent
In model theory, one defines the Grothendieck semigroup $\textup{Def(\S)} /_{\approx}$ of $\S$, where $\df{\S}:=\bigcup_n \Def{\S}$, and $a\approx b$ iff there is a definable bijection between $a$ and $b$, and where addition is induced by disjoint union.  In analogy with this, we call $\M/_{\simeq}$ 
the {\em Grothendieck semigroup of $\S$ with respect to $G$\/}, and $\M/_{\simeq_0}$
the {\em approximate Grothendieck semigroup (aGS) of $\S$ with respect to $G$\/} (with addition again induced by disjoint union in both cases).  We shall show that $aGS$ satisfies wCL: if $k [A] = k [B]$, then $[A] = [B]$.  
Stated differently,
$\M/_{\simeq}$, satisfies an approximate wCL, namely $k [A] \simeq_0 k [B]$ implies $[A] \simeq_0 [B]$.

In \cite{wagon}, p. 177, wCL in type semigroups is proved in the general (non-definable) setting. 
The proof uses the Hall-Rado-Hall infinite Marriage Theorem and is a variation on K\"{o}nig's method.  Essentially the same proof goes through here, except that we use Corollary \ref{multigraph} instead of the unavailable infinite Marriage theorem.
Here are the details.

\begin{theorem}[wCL]\label{cancellation}
Let $A,B \in \M$ be such that in $\M/_{\simeq}$, $k[A] \leq_0 k[B]$. Then $[A] \leq_0 [B]$.

\end{theorem}
\begin{proof}
We first assume that $A$, $B$ are 1-cycles.
Let $m>0$. Under the given hypothesis, we wish to find a definable $Y_{A} \subseteq A$ such that $Y_{A}\leq_m A$, and $A\setminus Y_{A} \leq B$.

Let \[kA = A_1 + A_2 + \dots + A_k,\] where $A_1=A$ and, for each $i\in [k]$, $A_1 \simeq A_i$ is witnessed by the pw-isomorphism $\phi_i$, with $\phi_1 = \textup{Id}_{A_1}$. Similarly, \[kB = B_1 + B_2 + \dots + B_k ,\] where $B_1 =B$ and, for each $i\in [k]$, $B_1 \simeq B_i$ is witnessed by the pw-isomorphism $\psi_i$, where $\psi_1 = \textup{Id}_{B_1}$.

Take $l>mk$ and a pw-embedding $\theta$ witnessing $kA \leq kB$ outside of a cycle $\Sigma_{i=1}^{k} V_i$ such that $\Sigma_{i=1}^{k} V_i \leq_l kA$ and where $V_i \subseteq A_i$, for each $i\in [k]$. 

Let $\mathcal{H}$ be the mp multigraph with bipartition $A , B$
and edge $(a,b)$ iff there are $1\leq i,j \leq k$ such that $\theta \circ \phi_i (a)=\psi_j (b)$.
Then $\mathcal{H}$ is of degree $k$ in $A$ outside of a definable subset $V$ such that $V\leq_{l} A$, and of degree $\leq k$  in $B$.  By Corollary \ref{multigraph}, there is a definable matching $M$ covering $A$ outside of a definable $Y_{A} \subseteq A$, where   $Y_{A } \leq_{m} A$.  Note that each of potentially multiple edges between $a\in A$ and $b\in B$ arose due to $\theta \circ \phi_i (a) = \psi_j (b)$ for a unique pair $(i,j)$, where $i,j \in [k]$.
Let \[C_{ij} = \{ a\in A \setminus Y_{0} \colon \, \exists b \big( b\in B \,\& \, (a,b)\in M \mbox{ due to } \theta \circ \phi_i (a) = \psi_j (b) \big) \} .\]
Then $\{ C_{ij} \}$ partitions $A\setminus Y_A$. The  
conclusion of the theorem is witnessed by the partition $\{ C_{ij}\}\cup Y_A$, and the maps $\psi_{j}^{-1} \circ \theta \circ \phi_i \colon C_{ij} \to B$.  Clearly, this map is injective.

Now assume $A=\Sigma_{i\in [m]}A_i$, $B=\Sigma_{j\in [n]}B_j$. To reduce to the case when $A$, $B$ are 1-cycles, we replace $A$ by $A_1 \times \{ 1\} \cup \dots \cup A_m \times \{ m \}$ and $B$ by $B_1 \times \{ 1 \} \cup \dots \cup  B_n \times \{ n \}$.  It is easy to see that the maps $\phi_i$, $\psi_i$ and $\theta$ can be modified accordingly so that the above proof goes through.
\end{proof}
Theorem \ref{aBSB} immediately yields:
\begin{corollary}\label{1sidedaWCL}
Let $[A] , [B] \in \M/_{\simeq}$ be such that $k[A] \simeq_0 k [B]$.  Then $[A] \simeq_0 [B]$.
\end{corollary}
We shall call a definable set $X$ {\em definably paradoxical\/} if $2X \simeq_0 X$ (this terminology will be justified by Theorem \ref{apprtarski}).

Here is a classical use of BSB and wCL:
\begin{proposition}
Suppose $X$ is definably paradoxical, and $Y$ is a definable subset of $X$ such that for some $k>0$, $X\leq_0 kY$.  Then $Y$ is also definably paradoxical.
\end{proposition}
\begin{proof}
We have 
\[ kX \simeq_0 X \leq_0 kY, \]
where the first relation follows from 
$2X\simeq_0 X$. Since also $kY \leq_0 kX$, we have $kX \simeq_0 kY$ by BSB. Then WCL yields $X\simeq_0 Y$, hence $Y$ is definably paradoxical. 
\end{proof}

\end{section}

\begin{section}{Equidecomposability}
In \cite{charges}, Laczkovich proves a criterion for equidecomposability in terms of charges (i.e. finitely additive measures that are allowed to take negative values), in the presence of the WCL.  The setting is that of a space $(X,G,\mathcal{A})$, where $X$ is a non-empty set, $G$ is a group of bijections $X\to X$, and $\mathcal{A}$ is a $G$-invariant boolean algebra on $X$ (i.e. a boolean algebra of subsets of $X$ such that if $a\in \mathcal{A}$ and $g\in G$, then $ga\in \mathcal{A}$).

\begin{definition}
A {\em $G$-invariant charge\/} is a map $\nu \colon \mathcal{A}\to \mathbb{R}\cup \{ \infty \}$ such that $\nu(gA)=\nu (A)$ for every $g\in G$, and $\nu(A \dot\cup B)=\nu (A)+\nu (B)$ for every $A,B \in \mathcal{A}$, where $\infty + \infty = \infty +a$ for every $a \in \mathbb{R}$.
\end{definition}
The main result in \cite{charges} is:
\begin{theorem}[\cite{charges}, p. 423, Corollary.]
In every space $(X,G,\mathcal{A})$, Conditions (1) and (2) below are equivalent iff the wCL holds (i.e. if $ka = kb$ implies $a=b$ for all $a,b \in \mathcal{A}$ and $k>0$):
\begin{enumerate}
\item Two elements $a,b$ from $\mathcal{A}$ are equidecomposable with respect to $G$.
\item For every $G$-invariant charge $\phi$ on $\mathcal{A}$, $\phi (a ) = \phi (b)$.
\end{enumerate}
\end{theorem}
That (1) implies (2) is obvious.  The proof of the non-trivial direction uses the following theorem.
\begin{theorem}[\cite{charges}, p. 423, Theorem 3.]\label{L2}
Let $(S,+)$ be a commutative semigroup, and let $a,b \in S$ be such that $\phi (a) =\phi (b)$ for every homomorphism $\phi \colon S \to \mathbb{R} \cup \{ \infty \}$. Then there is a positive integer $n$ such that $na=nb$.
\end{theorem}
This yields in our setting:
\begin{theorem}
Let $A,B \in \M$ be 1-cycles. Then $A \simeq_0 B$ iff $\phi A = \phi B$ for every $G$-invariant charge $\phi$ on $\Def{M}$.
\end{theorem}
\begin{proof}
The left-to-right direction is again obvious. For the other direction, 
take $\M/_{\simeq}$ with structure induced from $\M$ for $(S,+)$ in Theorem \ref{L2}, and use Corollary \ref{1sidedaWCL}.
\end{proof}
\end{section}

\begin{section}{Tarski's theorem}\label{tt}
Theorem \ref{cancellation} yields a quick proof of a variant of Theorem \ref{tarski}. 
To see this, we need a corollary of the wCL (see e.g. \cite{wagon}, p. 175, Cor. 10.21), whose proof we include here for the reader's convenience:
\begin{lemma}\label{from_k_to_2}
Let $[X]\in \M/_{\simeq}$, and suppose $(n+1)[X] \leq_0 n[X]$.  
Then $2[X] \leq_0 [X]$. 
\end{lemma}
\begin{proof}
By substituting $(n+1)X \leq_0 nX$ into itself, we obtain
\[ nX \geq_0 (n+1)X = nX + X \geq_0 (n+1)X + X = nX + 2X,\]
and after repeating this finitely many times,
\[ nX \geq_0 nX + nX = 2(nX). \]
Now apply Corollary \ref{1sidedaWCL} to $n(2X) \leq_0 nX$ to obtain $2X\leq_0 X$. 
\end{proof}

\smallskip\noindent
\begin{theorem}\label{apprtarski}
For a 1-cycle $X \in \M$, the following are equivalent.
\begin{enumerate}
\item $X$ is not definably paradoxical.
\item There is a $G$-invariant finitely additive measure $\mu \colon \M \to [0,\infty ]$ with $\mu X =1$.
\end{enumerate}
\end{theorem}
\begin{proof}
That (2) implies (1) is clear.  To prove the other implication, by Theorem \ref{tarski}, it suffices to show that $2X \not\leq_0 X$ implies $(k+1)X \not\leq kX$ for every $k\geq 1$.  Equivalently, that
$(k+1)X \leq kX$ for some $k\geq 1$ implies $2X \leq_0 X$.  But that is just Lemma \ref{from_k_to_2}.
\end{proof}
Finally, we would like to mention some related work.  In \cite{k},
Kraj\'{i}\v{c}ek denotes by $WPHP^{2n}_{n}$ the condition that there is no definable set such that two copies of it embed into one, using definable bijections, and characterizes structures satisfying $WPHP^{2n}_{n}$ in terms of the existence of an approximate Euler characteristic (aEC), which is, roughly, a map from $\df{\S}$ into a partially ordered ring, in which the usual equalities for an Euler characteristic are replaced by approximate versions, and which factors through the quotient map $\df{\S} \to \df{\S}/_{\approx}$. 
In particular, the paper proves that if $\S$ satisfies $WPHP^{2n}_{n}$ and $CC$ (for every $a,b \in \Def{\S}$, $a$ embeds into $b$, or vice versa), then $\S$ admits a non-trivial approximate Euler characteristic, where the target ring is graded by a dimension function, and on each level the aEc is a Keisler measure (see Claim 1 in proof of Proposition 5.4, p. 575), invariant under definable bijections. 

Condition (1) in Theorem \ref{apprtarski} is also related to the notion of a $BpD$ in Chernikov et al. \cite{small} (which really goes back to \cite{nip}, where it appears in the group setting): Suppose a group $G$ acts on a set $S$, and $B$ is a $G$-invariant boolean algebra of subsets of $S$.  Then $S$ has $BpD$ if there is a pw isomorphic embedding (with pieces from $B$) $S+Y \to Y$ for some cycle $Y$. 
The following result from \cite{small} was already implicit in \cite{nip} (Fact 4.4) and is a special case of Theorem \ref{apprtarski} (where the distinguished set $X$ is $S$):
\begin{proposition}[\cite{small}, Prop 4.11]\label{BpD} There is a $G$-invariant finitely additive probability measure on $B$ iff $S$ has no $BpD$.  
\end{proposition}

\end{section}

\end{document}